\documentclass[10pt, reqno]{amsart}
\usepackage{amsmath, amsthm, amscd, amsfonts, amssymb, graphicx, color}
\usepackage[bookmarksnumbered, colorlinks, plainpages]{hyperref}
\hypersetup{colorlinks=true,linkcolor=red, anchorcolor=green, citecolor=cyan, urlcolor=red, filecolor=magenta, pdftoolbar=true}
\usepackage{mathrsfs}

\theoremstyle{definition}

\theoremstyle{remark}

\begin{document}
\setcounter{page}{1}

\title[Comments on a note by M.~Waldschmidt]{Comments on a note by M.~Waldschmidt}

\author[Nikolaev]
{Igor V. Nikolaev$^1$}

\address{$^{1}$ Department of Mathematics and Computer Science, St.~John's University, 8000 Utopia Parkway,  
New York,  NY 11439, United States.}
\email{\textcolor[rgb]{0.00,0.00,0.84}{igor.v.nikolaev@gmail.com}}





\begin{abstract}
This  note is the follow up to a paper by M.~Waldschmidt. 
\end{abstract}

\maketitle

In \cite{Wal1} it was argued that the value of $\mu=\log \varepsilon$ in  \cite[Theorem 1]{Nik1}   must be equal to $1$. 
As a consequence, \cite[Theorem 1]{Nik1}  is claimed to be false.
In this note we explain why $\mu=\log \varepsilon$ 
cannot be  $1$,   thus upholding   \cite[Theorem 1]{Nik1}. 
Our comments touch  \cite[Second paragraph]{Wal1},  since   \cite[First paragraph]{Wal1} is covered in 
\cite[Remark 2]{Nik1}.

The $\beta$-KMS (Kubo–Martin–Schwinger) state on a $C^*$-algebra $\mathcal{A}$ is a real number $\beta$,
such that
\begin{equation}
\omega(x\sigma_{i\beta}(y))=\omega (yx), \quad\forall x,y\in\mathcal{A}, 
\end{equation}
where $\sigma_t$ is the 1-parameter group of  automorphisms of $\mathcal{A}$ and 
$\omega$ is a weight on $\mathcal{A}$. 
Such a state is 
 the fundamental invariant of the algebra $\mathcal{A}$
  known  in the quantum statistical mechanics as an inverse
temperature.
If $\mathcal{A}$ is an AF (Approximate Finite) $C^*$-algebra given by  a
constant incidence matrix $A\in GL_n(\mathbf{Z})$,  the inverse temperature
was calculated in  \cite[Chapter 4 and Fig. ~1-5]{BJO} 
by the remarkable formula
\begin{equation}\label{eq2}
\beta=\log ~(\lambda_A),
\end{equation}
where $\lambda_A>1$ is the Perron-Frobenius eigenvalue of the matrix $A$.  
Moreover,  the KMS state (\ref{eq2}) is unique, {\it ibid.}

In our case $n=2$ and the Pimsner-Voiculescu theorem says that the non-commutative torus
$\mathcal{A}_{\theta}$ embeds into an AF-algebra preserving the respective $K_0$-groups. 
In particular, the incidence matrix $A$ is nothing but the matrix form of the fundamental unit
$\epsilon>1$ of the number field $\mathbb{Q}(\theta)$ and therefore the inverse temperature of $\mathcal{A}_{\theta}$
\begin{equation}\label{eq}
\beta=\log ~(\epsilon). 
\end{equation}
Since  $\mu$  in \cite[Lemma 3]{Nik1} is known to parametrize  the state space of the 
$C^*$-algebra  $\mathcal{A}_{\theta}$,  we identify it with $\beta$ keeping the notation
in below.

Assuming $\beta=1$  \cite[Second paragraph]{Wal1}  would simply mean  that a representation of 
the Sklyanin $\ast$-algebra $\left\{S(q) ~|~ q=\beta e^{2\pi i\theta}\right\}$ by the linear operators on a Hilbert space does not exist,
 since there is no scalar product to apply the standard GNS (Gelfand–Naimark–Segal)  construction. 
 The same is  true for all other  values of $\beta$,  unless  $\beta=\log ~(\epsilon)$ in which case the GNS construction
 gives a representation of the Sklyanin algebra.  Such a representation is critical for the validity of \cite[Lemma 2]{Nik1}.

\bibliographystyle{amsplain}

\begin{thebibliography}{99}

\bibitem{BJO}
O.~Bratteli,  P. ~E. ~T. ~Jorgensen and  V. Ostrovskyı,
\textit{Representation Theory and Numerical AF-invariants.
The representations and centralizers of certain states on $\mathcal{O}_d$},
Memoirs of the A. M. S. {\bf 168} (2004),  No. 797 

\bibitem{Nik1}
I.~Nikolaev \textit{On algebraic values of function $\exp~(2\pi ix +\log\log y)$},  Ramanujan J. {\bf 47} (2018), 417–425.


\bibitem{Wal1}
M.~ Waldschmidt, \textit{On a paper by Nikolaev}, 
Ramanujan J. {\bf 57} (2022), 1517-1518.

\end{thebibliography}


\end{document}